\documentstyle[twoside,draft,amssymb,amsmath,amsfonts]{article}

\def\RR{\mathbb{R}}
\def\NN{\mathbb{N}}
\def\CC{\mathbb{C}}
\def\ZZ{\mathbb{Z}}
\def\PP{\mathbb{P}}

\newcommand{\eop}{\vrule height7pt width7pt depth0pt%

    \vrule height0pt width0pt depth20pt}

\begin{document}


\renewcommand{\maketitle}{

   \vskip 2em%

\begin{flushleft}\tiny

Trends and Applications in Constructive Approximation\\

(Eds.) M.G. de {\bf Bruin}, D.H. {\bf Mache} \& J. {\bf Szabados}\\

\smallskip

International Series of Numerical Mathematics Vol. 1??\\

\copyright 2005 Birkh\"auser Verlag Basel (ISBN 3-7643-7124-2)\\[2cm]

\end{flushleft}

  \begin{center}%

  \let \footnote \thanks

    {\LARGE \title \par}%

    \vskip 1.5em%

    {\large

      \lineskip .5em%

      \begin{tabular}[t]{c}%

        \author

      \end{tabular}\par}%

    \vskip 1em%

  \end{center}%

  \par

  \vskip 1.5em

  \setcounter{section}{0}

  \setcounter{subsection}{0}

  \setcounter{equation}{0}

  \thispagestyle{plain}}




\newcommand{\reff}[1]{(\ref{#1})}

      \setcounter{section}{0}

       \def\title{Recent progress on univariate and multivariate\\polynomial and spline quasi-interpolants}

						\def\author{Paul Sablonni\`ere}


                   {Polynomial and spline quasi-interpolants}

       \date{}

       \maketitle

       \begin{abstract}

       \noindent

      Polynomial and spline quasi-interpolants (QIs) are practical and
effective approximation operators.
Among their remarkable properties, let us cite for example: good shape
properties, easy
computation and evaluation (no linear system to solve), uniform boundedness
independently of the
degree (polynomials) or of the partition (splines), good approximation order.
We shall emphasize new results on various types of univariate and multivariate
polynomial or spline QIs, depending on the nature of coefficient
functionals, which can be
differential, discrete or integral. We shall also present some applications
of QIs to
numerical methods.

 \end{abstract}


\section{Introduction}

\noindent
A quasi-interpolant of $f$ has the general form
$$
Qf=\sum_{\alpha\in A} \mu_{\alpha}(f) B_{\alpha},
$$
where $\{B_{\alpha},\alpha\in A\}$ is a family of polynomials or B-splines
forming a partition of unity,
and $\{\mu_{\alpha}(f),\alpha\in A\}$
is a family of linear functionals which are local in the sense that they
only use values of $f$
in some neighbourhood of  $\Sigma_{\alpha}=supp(B_{\alpha})$. The main
interest of QIs is that they
provide excellent approximants of functions without solving any linear
system of equations.
In the literature, one can find the three following types of QIs:

(i) Differential QIs (abbr. DQIs) : the linear functionals are linear combinations of values of
derivatives of $f$ at some point in $\Sigma_{\alpha}$.

(ii) Discrete QIs (abbr. dQIs) : the linear functionals are linear
combinations of values of $f$
at some points in the neighbourhood of $\Sigma_{\alpha}$.

(iii) Integral QIs (abbr. iQIs) : the linear functionals are linear
combinations of
weighted mean values of $f$ in the neighbourhood of $\Sigma_{\alpha}$.

\noindent
We shall present various types of univariate and multivariate polynomial
and spline QIs, mainly
dQIs and iQIs, which were recently introduced in the literature. For
polynomial QIs, we only present
QIs which are close to the orifginal Bernstein or Durrmeyer operators (for
other types of QIs,
see for example \cite{GM2}\cite{GM3}).

The prototype of {\sl polynomial dQIs} is the classical Bernstein operator
$$
B_nf=\sum_{i=0}^{n} f(\frac{i}{n}) b_i^{(n)}
$$
where $\{b_i^{(n)}(x)=C_n^i x^i(1-x)^{n-i}, 0\le i\le n\}$ is the Bernstein
basis of the space
${\PP}_n$ of polynomials of degree at most $n$ (the $C_n^i$ are binomial
coefficients).

The prototype of {\sl polynomial iQIs} is the Durrmeyer operator \cite{Dur}
$$
M_nf=\sum_{i=0}^{n} \langle f,\tilde b_i^{(n)} \rangle b_i^{(n)}
$$
where $\tilde b_i^{(n)}=b_i^{(n)}/\int_0^1 b_i^{(n)}=(n+1) b_i^{(n)}$ and
$\langle f, g \rangle=\int_0^1 fg$.
Both can be extended to the multivariate case, either on the hypercube or
on the simplex.
Another extension consists in adding a Jacobi weight in the scalar product.

The prototypes of {\sl spline DQIs} are de Boor-Fix QIs \cite{dBF} and
their various univariate and multivariate
extensions
$$
Qf=\sum_{j\in J} \lambda_j(f) B_j.
$$
Here $\{B_j\; j\in J\}$ is a family of univariate B-splines of degree $m$
on a nonuniform
sequence of knots $\{t_k\}$. Assuming that
$\Sigma_j=supp(B_j)=[t_{j-m},t_{j+1}]$, we set
$E_m=\{-m+1,\ldots,0\}$ and we define
$\psi_j(t)=\prod_{r\in E_m}(t_{j+r}-t)\in{\PP}_{m}$ for all
$j\in J$. For any $\tau\in \Sigma_j$, the coefficient functionals are
$$
\lambda_j(f)=\frac{1}{(m-1)!}\sum_{l=0}^{m-1}
(-1)^{m-l-1}D^{m-l-1}\psi_j(\tau)D^lf(\tau).
$$
The prototypes of {\sl spline dQIs} are the various univariate and
multivariate extensions of
Schoenberg-Marsden operators \cite{MS1}\cite{MS2}.
$$
Sf=\sum_{j\in J} f(\tau_j) B_j
$$
where  $\tau_j$ is an interior point of $\Sigma_j=supp(B_{\alpha})$.

The prototypes of {\sl spline iQIs} are the various univariate and
multivariate extensions
of operators \cite{Cz}\cite{Sab3}
$$
Tf=\sum_{j\in J} \langle f, M_j \rangle B_j,
$$
where $M_j$ is a B-spline (which can be different from $B_j$) normalized by
$\int_{\Sigma_j}  M_j=1$.\\

\noindent
As emphasized by de Boor (\cite{dB3}, chapter XII), a spline QI defined on
non uniform partitions
has to be {\sl uniformly bounded independently of the partition} (abbr. UB)
in order to be
interesting for applications. Therefore, with some coworkers, we have
defined various families of
QIs satisfying this property and having  an infinite norm as  small as
possible. In general it is
difficult to minimize the true norm of the operator,  however, it is often
possible to minimize
an upper bound of this norm: this gives rise to what we have called
near-best (abbr. NB) QIs
(see \cite{BIS},\cite{BISS1}-\cite{BISS3},\cite{IP}).\\
Numerical applications are still not very much developed. However, QIs can
be useful in
approximation and estimation \cite{Cz}\cite{CGS}\cite{Wal}, in numerical
quadrature
\cite{DL}\cite{Sab13}\cite{Sab15}, and for the numerical solution of
integral or partial
differential equations.

----------------------------------

\section{Univariate polynomial QIs }

\subsection{Basic operators}

1) The {\sl Bernstein-Stancu QI} \cite{Sta} is defined for $x\in [0,1]$ by
$$
S_n^{(\alpha)}f(x)=\sum_{i=0}^{n} f(\frac{i}{n}) b_i^{(n)}(x,\alpha)
$$
where the {\sl Bernstein-Stancu basis} is defined by
$\displaystyle
b_k^{(n)}(x,\alpha)=C_n^k\frac{(x)_{\alpha}^k(1-x)_{\alpha}^{n-k}}{(1)_{\alpha}^n}.$

Here $(x)_{\alpha}^k=x(x+\alpha)\ldots (x+(k-1)\alpha)$,for $\alpha\in \RR$.
For $\alpha=0$, we recover the classical Bernstein basis.

2) The {\sl Bernstein- Phillips} (or {\sl q-Bernstein) QI}
(\cite{OP}-\cite{Ph4}) is
defined for $x\in [0,1]$ by
$$
B_n^{q}f(x)=\sum_{i=0}^{n} f\left(\frac{[i]}{[n]}\right) b_k^{(n)}(x,q)
$$
where the {\sl Bernstein-Phillips or q-Bernstein basis} is defined for
$q\ne 1$ by\\
$b_k^{(n)}(x,q)=\Gamma_n^k x^k (1-x)_q^{n-k}$.
Here $[i]=\frac{1-q^i}{1-q}$, $[i]!=\prod_{s=1}^{i} [s]$,
$\Gamma_n^k=\frac{[n]!}{[k]![n-k]!}$ and $(x)_q^k=\prod_{s=0}^{k-1}(1-q^sx)$.
For $q=1$, we recover the classical Bernstein basis.\\

Using the notation $e_s(x)=x^s$ for monomials, it is easy to prove that all
the above QIs
${\mathcal{B}}_n$ are exact on ${\PP}_1$, i.e. ${\mathcal{B}}_n e_s=e_s
\;\; for \;\; s=0,1$.
Moreover they are degree preserving since ${\mathcal{B}}_n
e_s(x)=e_s(x)+r_{s-1}(x,n)$
where $r_{s-1}(x,n)$ is some polynomial of degree at most $s-1$ depending
on $n$ (and
eventually on the parameters $\alpha$, $N$ or $q$).

QIs------------------------------------

\subsection{Left and right BQIs}

All operators ${\mathcal{B}}_n$ defined above are isomorphisms of
${\PP}_n$. Moreover
${\mathcal{B}}_n$ and ${\mathcal{A}}_n={\mathcal{B}}_n^{-1}$  can be
expressed as
linear differential operators with polynomial coefficients
$$
{\mathcal{B}}_n=\sum_{k=0}^{n} \beta_k^{(n)}D^k,\quad
{\mathcal{A}}_n=\sum_{k=0}^{n}
\alpha_k^{(n)}D^k,
$$
where $D=\frac{d}{dx}$ and the polynomials $\beta_k^{(n)}\in {\PP}_k$ and
$\alpha_k^{(n)}\in {\PP}_k$
are defined by simple recursions (see e.g.\cite{Sab4}-\cite{Sab7} for
partial results in this
sense).

For $0\le r\le n$, we introduce the partial inverses:
$$
{\mathcal{A}}_n^{(r)}=\sum_{k=0}^{r} \alpha_k^{(n)}D^k,
$$
and we consider the two families of right and left BQIs:\\

\noindent
(RBQI) The right BQIs
${\mathcal{B}}_n^{[r]}={\mathcal{B}}_n\circ {\mathcal{A}}_n^{(r)}$
are defined for $C^r$-functions $f$ by
$$
{\mathcal{B}}_n^{[r]}f={\mathcal{B}}_n ( {\mathcal{A}}_n^{(r)}
f)={\mathcal{B}}_n (\sum_{k=0}^{r}
\alpha_k^{(n)} D^k f).
$$
\noindent
(LBQI) The left BQIs
${\mathcal{B}}_n^{(r)}={\mathcal{A}}_n^{(r)}\circ {\mathcal{B}}_n$
are defined on any (e.g. continuous) function
$$
{\mathcal{B}}_n^{(r)}f={\mathcal{A}}_n^{(r)}({\mathcal{B}}_n f)=\sum_{k=0}^{r}
\alpha_k^{(n)}D^k({\mathcal{B}}_n f).
$$

By construction, for $0\le r\le n$, the BQIs ${\mathcal{B}}_n^{[r]}$ and
${\mathcal{B}}_n^{(r)}$ are
exact on the space ${\PP}_r$. Moreover, in many cases, the LBQIs have a
uniformly bounded infinite
norm, independent on $n$ for each $0\le k\le  n$ fixed (see e.g.
\cite{Sab10} \cite{WZ} for
some results of this type). From this property are deduced some convergence
results (see
\cite{Dia}\cite{Sab7}).

-------------------------------------

\subsection{Kageyama QIs}

Kageyama \cite{Kag1} \cite{Kag2} considers Stancu operators for $\alpha\in
[-\frac1n,0]$
$$
S_n^{(-\frac1n)}={\mathcal{L}}_n, \;\; S_n^{(0)}=B_n,\;\;
$$
where ${\mathcal{L}}_n$ is the Lagrange interpolation operator on the
uniform partition of $[0,1]$
(this result is due to M\"uhlbach). Then he truncates at order $s$ the
Maclaurin series of
$S_n^{(\alpha)}f$ w.r.t. $\alpha$ and  he takes the value of this
polynomial at $\alpha=-\frac1n$:
$$
{\mathcal{K}}_n^{(s)}f=\sum_{j=0}^{s}\frac{1}{j!}\frac{(-1)^j}{n^j}
\frac{\partial^j}{\partial\alpha^j}\left[ S_n^{(\alpha)}f \right]_{\alpha=0}
$$
${\mathcal{K}}_n^{(0)}=S_n^{(0)}=B_n$ and
${\mathcal{K}}_n^{(\infty)}=S_n^{(-\frac1n)}={\mathcal{L}}_n.$
He also gives expansions of ${\mathcal{K}}_n f$ in terms of derivatives
of $B_nf$ and in powers of $\frac1n$.
He proves that, for all $s$ fixed, $\Vert {\mathcal{K}}_n^{(s)}
\Vert_{\infty}$ is uniformly bounded
and give Voronovskaja type results, e.g.
$$
\lim_{n\to\infty}n^{s+1}({\mathcal{K}}_n^{(s)}f-f)=
-\sum_{k=0}^{2s+2}\frac{1}{k!}\Upsilon_{s+1,k}D^kf
$$
where the polynomials $\Upsilon_{s+1,k}$ can be computed by recursion.
He also compares the expansions of ${\mathcal{L}}_n$, the BQIs $B_n^{(r)}$ and ${\mathcal{K}}_n^{(s)}$
interms of derivatives of $B_nf$ with polynomial coefficients. Numerical
experiments done by
the author suggest that  these operators are in general better approximants
than BQIs of section
3.1.

------------------

\subsection{Univariate Durrmeyer and Goodman-Sharma QIs}

A straightforward generalization of the Durrmeyer operator $M_n$ consists
in introducing a Jacobi weight on $[0,1]$
in the associated scalar product
$$
\langle f,g\rangle=\int_0^1 w_{\alpha,\beta}(t)f(t)g(t)dt,\;\;
w_{\alpha,\beta}(t)=t^{\alpha}(1-t)^{\beta}, \;\; for \;\; \alpha,\beta>-1
$$
The extended Durrmeyer-Jacobi operator (\cite{BX},\cite{Sab1}) is then
defined by
$$
M_n^{(\alpha,\beta)}f=\sum_{i=0}^{n} \frac{\langle f,b_i^{(n)}
\rangle}{\langle e_0, b_i^{(n)}.
\rangle} b_i^{(n)}
$$
The limit case $(\alpha,\beta)=(-1,-1)$, corresponding to the weight
$\displaystyle \tilde w(x)=\frac{1}{x(1-x)}$, gives a QI with very
attractive properties.
It has been introduced by Goodman and Sharma \cite{GS1}\cite{GS2} for
polynomial (and a variant
for  spline) QIs. It can be written as follows, with $Lf(x)=(1-x)f(0)+xf(1)$:
$$
G_nf=Lf+(n-1)\sum_{i=1}^{n-1}\langle f-Lf,b_{i-1}^{(n-2)}\rangle b_i^{(n)}.
$$
This operator is exact on ${\PP}_1$ and its behaviour is quite similar to
that of the classical Bernstein operator.
For example, one has for $f\in C^2(I)$
$$
\lim n(f(x)-G_nf(x))=x(1-x)f''(x).
$$
It also preserves the positivity, the monotonicity and the convexity of
$f$. As discrete
Bernstein operators, the above  operators $G_n$ have associated QIs in the
sense of section
3.1 \cite{Sab10}
.

\subsection{Extrapolation}

All operators ${\mathcal{B}}_n$ described in this section have asymptotic
expansions of type
$$
{\mathcal{B}}_nf(x)\sim f(x)+\sum_{k\ge r}\frac{\varphi_k^{(n)}(f,x)}{n^k}
$$
for some index $r$, the $\varphi_k^{(n)}(f,x)$ being linear differential
operators depending on $n$ and $k$.
Therefore they are good candidates for extrapolation methods (see e.g.
\cite{BRZ} and
\cite{Sid}).  Numerical experiments done by the author show that Richardson
extrapolation is
efficient while the use of variants of epsilon or
$\Delta^2$ algorithms often introduce spurious poles in the interval of
definition.

--------------------------

\section{Polynomial QIs on a simplex}

\subsection{Bernstein operator and associated QIs}

The simplex $S$ of dimension $d-1$ is defined in barycentric coordinates as

\noindent
$S=\{{\mathbf{x}}=(x_1,x_2,\ldots,x_d): \vert {\mathbf{x}} \vert=1\}$
with $\vert {\mathbf{x}} \vert=\sum_{i=1}^d \vert x_i\vert$.

\noindent
The associated simplex of indices, monomials and partial derivatives are
defined by
$\Sigma_n=\{{\mathbf{i}}=(i_1,i_2,\ldots,i_d):
\vert {\mathbf{i}}\vert=n\}$,
$X_n=\{\frac{\mathbf{i}}{n}: {\mathbf{i}}\in \Sigma_n\}\subset S$,

\noindent
${\mathbf{i}}!=i_1!i_2!\ldots,i_d!,\;\;\;\;$
${\mathbf{x}}^{\mathbf{i}}=x_1^{i_1}x_2^{i_2}\ldots x_d^{i_d}$,
$D^{\mathbf{i}}=D_1^{i_1}D_2^{i_2}\ldots D_d^{i_d}$
with $D_s=\frac{\partial}{\partial x_s}.$

\noindent
The Bernstein basis of ${\PP}_n$ (space of polynomials of total degree at
most $n$) and
the Bernstein operator are defined respectively by:
$$
b_{\mathbf{i}}^{(n)}({\mathbf{x}})=\frac{n!}{\mathbf{\mathbf{i}}!}{\mathbf{x
}}^{\mathbf{i}}\quad
for \;\;{\mathbf{i}}\in \Sigma_n,\;\;
B_n f({\mathbf{x}})=\sum_{{\mathbf{i}}\in
\Sigma_n} f(\frac{{\mathbf{i}}}{n}) b_{{\mathbf{i}}}^{(n)}(\mathbf{x})
$$
As $\sum_{{\mathbf{i}}\in\Sigma_n} b_{\mathbf{i}}^{(n)}({\mathbf{x}})=1$ and
$\sum_{{\mathbf{i}}\in \Sigma_n}\frac{i_s}{n}
b_{\mathbf{i}}^{(n)}({\mathbf{x}})={\mathbf{x}}^{\varepsilon_s}=x_s$, for
$1\le s\le d$,
where\\ $\varepsilon_s=(0,0,\ldots,1,\ldots,0)$, then
$B_n$ is exact on ${\PP}_1$.

Let  ${l_{\mathbf{i}}^{n},\mathbf{i}\in \Sigma_n}\}$
be the Lagrange basis of ${\PP}_n$ associated with the data
points $X_n$. Then
$l_{\mathbf{i}}^{n}(\frac{\mathbf{j}}{n})=\delta_{\mathbf{i}\mathbf{j}}$
implies
$B_n l_{\mathbf{i}}^{n}=b_{\mathbf{i}}^{(n)}$, hence
$B_n$ is an isomorphism of ${\PP}_n$.

\noindent
For $f\in C^2(S)$, we have the Voronovskaja type result
(\cite{Lai2}\cite{Sau}\cite{Wal}),
$$
\lim n\left[B_nf-f\right]=\frac12 \bar{\mathcal{D}}f
$$
where ${\mathcal{D}} f$ is the differential operator
$$
\bar{\mathcal{D}} f(x)=\sum_{i<j} x_ix_j (\partial_i-\partial_j)^2.
$$
$B_n$ and  its inverse $A_n=B_n^{-1}$ in ${\PP}_n$ can be expressed as
linear differential operators
$$
B_n=\sum_{\mathbf{i}\in \Sigma_n} \beta_{\mathbf{i}}^{(n)}D^{\mathbf{i}},\quad
A_n=\sum_{\mathbf{i}\in \Sigma_n} \alpha_{\mathbf{i}}^{(n)}D^{\mathbf{i}}
$$
whose coefficients can be computed by recursion.
For $0\le k\le n$, define partial inverses
$$
A_n^{(k)}=\sum_{\mathbf{i}\in
\Sigma_k}\alpha_{\mathbf{i}}^{(n)}D^{\mathbf{i}}.
$$
As in section 2.2 for univariate QIs, we can consider the two families of
operators:

left Bernstein quasi-interpolants (LBQIs) $B_n^{(k)}=A_n^{(k)}\circ B_n$, and
right Bernstein quasi-interpolant (RBQIs) $B_n^{[k]}=B_n\circ A_n^{(k)},$
where $B_n^{(0)}=B_n^{[0]}=B_n=$ and $B_n^{(n)}=B_n^{[n]}={\mathcal{L}}_n=$
Lagrange  interpolation on $X_n$.

We have proved \cite{Sab5} that
$\Vert B_n^{(2)}\Vert_{\infty}\le 2d+1$ for all $n\ge 2$, and we
conjecture that for all $k\ge 0$, there exists a constant $C_k(d)$ such
that for all
$n\ge k$,
$$
\Vert B_n^{(k)}\Vert_{\infty}\le C_k(d).
$$
We also conjecture the Voronovskaja-type results
$$
\lim n^{r+1}(B_n^{(2r)}f-f)={\mathcal{A}}_{2r}f,\;\;
\lim n^{r+1}(B_n^{(2r+1)}f-f)={\mathcal{A}}_{2r+1}^*f,
$$
where ${\mathcal{A}}_{2r}$ and ${\mathcal{A}}_{2r+1}^*$ are linear
differential operators,
and the asymptotic expansions
$$
B_n^{(2r)}f\; and \; B_n^{(2r+1)}f \sim
f+\frac{c_{r+1}}{n^{r+1}}+\frac{c_{r+2}}{n^{r+2}}\ldots
$$
--------------------------

\subsection{Durrmeyer-Jacobi  QIs on a simplex}

One can introduce a Jacobi weight on the simplex in the scalar product of
$L_w^2(S)$:
$w_{\alpha}(x)=x^{\alpha},\;\; \langle
f,g\rangle=\int_Sw_{\alpha}(x)f(x)g(x)dx,$
and define the Durrmeyer-Jacobi quasi-interpolants (DJQIs)
$$
M_nf=\sum_{i\in \Sigma_n}\frac{\langle f,b_i^{(n)}\rangle}{\langle e_0,
b_i^{(n)}\rangle}b_i^{(n)}.
$$
Its eigenvectors are the Jacobi polynomials on the simplex.
There holds a \\Voronovskaja type result \cite{BrS}\cite{Sau}
$$
\lim n(M_nf(x)-f(x))={\mathcal{D}}_\alpha f(x)
$$
where the differential operator $\mathcal{D}_\alpha$ is defined by
$$
{\mathcal{D}}_\alpha=x^{-\alpha}\sum_{i<j} (\partial_i-\partial_j)
x_ix_j^{\alpha}
(\partial_i-\partial_j).
$$
As $M_n$ is an isomorphism of ${\PP}_n$, one can expand
$M_n=\sum_{k=0}^n\sum_{\mathbf{i}\in
\Sigma_k}\beta_{\mathbf{i}}^{(n)}D^{\mathbf{i}}$ and
 $L_n=M_n^{-1}=\sum_{k=0}^n\sum_{\mathbf{i}\in \Sigma_n}
\alpha_{\mathbf{i}}^{(n)}D^{\mathbf{i}}.$
As in the univariate case \cite{Sab11}, the polynomials
$\beta_{\mathbf{i}}^{(n)}$ and
$\alpha_{\mathbf{i}}^{(n)}$ are probably linear combinations of Jacobi
polynomials on $S$
(\cite{Der3})
Setting $L_n^{(r)}=\sum_{k=0}^r\sum_{\mathbf{i}\in \Sigma_n}
\alpha_{\mathbf{i}}^{(n)}D^{\mathbf{i}}$, one can define the left DJQIs
$M_n^{(r)}=L_n^{(r)}\circ M_n$,  and the right DJQIs $M_n^{[r]}=M_n\circ
L_n^{(r)}$, with
$M_n^{(0)}=M_n$ and $M_n^{(n)}=P_n=$ orthogonal projector on ${\PP}_n$ in
$L^2(S)$. They have the same properties as univariate QIs, and it would be
interesting to
have detailed proofs, those of \cite{Sab5}\cite{Sab6} being only sketched.
However, the author thinks that the following operators are still more
attractive.

-----------
\subsection{Jetter-St\"ockler operators on a triangle}

For the sake of simplicity, we describe them over a triangle (with
barycentric coordinates $\{\lambda_1,\lambda_2,\lambda_3\}$)
in the case of the Legendre weight ($w=1$, see \cite{JS2} for the general
study on a simplex with Jacobi weight).
Using the following notations:
$$
D_{ij}=\partial_j-\partial_i, \; i<j,\;\; D=\{D_{12},D_{13},D_{23}\},\;
\Lambda=\{\lambda_1\lambda_2,\lambda_1\lambda_3,\lambda_2\lambda_3\},\;
$$
$$
{\mathbf{k}}=(k_{12},k_{13},k_{23})\in{\NN}^3,\;
{\mathcal{D}}^{\mathbf{k}}=D_{12}^{k_{12}}D_{13}^{k_{13}}D_{23}^{k_{23}},
$$
$$
\Lambda^{\mathbf{k}}=(\lambda_1\lambda_2)^{k_{12}}(\lambda_1\lambda_3)^{k_{1
3}}(\lambda_2\lambda_3)^{k_{23}},
$$
the authors define the following basic differential operators:
$$
U_{\mathbf{k}}=\frac{1}{{\mathbf{k}}!}(-1)^{\vert\mathbf{k}
\vert}{\mathcal{D}}^{\mathbf{k}}\Lambda^{\mathbf{k}}{\mathcal{D}}^{\mathbf{k}},
\;\;
{\mathcal{U}}_\ell=\frac{1}{\ell!}\sum_{\vert\mathbf{k}
\vert=\ell}U_{\mathbf{k}},
\;\;
{\mathcal{Y}}_{n}=\sum_{\ell=0}^{n}(C_{n}^{\ell})^{-1}{\mathcal{U}}_\ell
$$
Let $M_n$ be the Durrmeyer operator, then they prove that
$U_{\mathbf{k}}$ commute with $M_n$ for all pairs $({\mathbf{k}},n)$ and that
${\mathcal{Y}}_{n}$ is
the inverse of $M_n$ in the space of polynomials ${\PP}_n$. Now,  for $0\le
r\le n$ fixed, they
define partial inverses and left Jetter-St\"ockler quasi-interpolants (LJSQIs)
$$
{\mathcal{Y}}_n^{(r)}=\sum_{\ell=0}^{r}(C_n^{\ell})^{-1}{\mathcal{U}}_\ell,\;\;
M_n^{(r)}={\mathcal{Y}}_n^{(r)}M_n.
$$
One can also define right JSQIs $M_n^{[r]}=M_n {\mathcal{Y}}_n^{(r)}$.
Both operators $M_n^{(r)}$ and $M_n^{[r]}$ are exact on ${\PP}_r$.
Moreover, for $r$ fixed, the left JSQIs have uniformly bounded infinite
norms w.r.t. $n$.
Finally, the authors  prove Voronovskaja-type results:
$$
\lim_{n\to\infty}C_n^r(f-M_n^{(r-1)}f)={\mathcal{U}}_r f
$$

\subsection{Extrapolation}

All operators ${\mathcal{B}}_n$ described in this section have asymptotic
expansions of type
$$
{\mathcal{B}}_nf(x)\approx f(x)+\sum_{k\ge r}\frac{\varphi_k^{(n)}(f,x)}{n^k}
$$
(see e.g. \cite{Lai2} and \cite{Wal}). In particular, the latter reports
interesting numerical
results on Richardson extrapolation of classical Bernstein operators on the
triangle.
It would be interesting to compare these results with those which could be
obtained by extrapolating the above QIs.

PARTITIONS-----------------------

\section{Univariate spline QIs on uniform partitions}

\subsection{Univariate differential and discrete QIs}

For the construction of QIs with optimal approximation order, we refer to
\cite{BSSV} and \cite{BS}, where  general solutions are given, thus completing
 the initial work by Schoenberg in \cite{Sch}.

\subsection{Near-best spline dQIs}

Consider the family of spline dQIs of order $2m$ depending on
$n+1$ arbitrary parameters $a=(a_0,a_1,\ldots,a_n)$, $n\ge m$:
$$
Q_af=\sum_{i\in\ZZ} \Lambda f(i) M_{2m}(x-i)
$$
with coefficient functionals
$$
\Lambda f(i)=a_0f(i)+\sum_{j=1}^{n}a_j\left(f(i+j)+f(i-j)\right).
$$
Setting $\nu(a)=\vert a_0 \vert+\sum_{j=1}^{n}\vert a_j\vert$,
then we have $\Vert Q_a\Vert_{\infty}\le \nu(a)$.
By imposing that $Q_a$ be exact on ${\PP}_r$, with $0\le r\le 2m-1$,
we obtain a set of linear constraints: $a\in V_r\subset {\RR}^{n+1}$.
We say that $Q^*=Q_{a^*}$ is a {\sl near best dQI} if
$$
\nu(a^*)=\min\{\nu(a); a\in V_r\}.
$$
There is existence, but in general not unicity, of solutions.

\noindent
{\sl Example:} cubic splines (see \cite{IP}).
There is a unique optimal solution for $n\ge 2$:
$$
a_0^*=1+\frac{1}{3n^2},\;\; a_n^*=-\frac{1}{6n^2},\;\; a_j^*=0 \;\; for
\;\; 1\le j\le n-1
$$
Moreover, for all $n\ge 4$, $\Vert Q^* \Vert_{\infty}\le 1+\frac{2}{3n^2}.$
Here are the first values of $\Vert Q^* \Vert_{\infty}$ and $\nu(a^*);$
$n=1: 1.222\;\; \&\;\; 1.666;$
$n=2: 1.139\;\;  \&\;\; 1.166;$
$n=3: 1.074\;\;  \&\;\; 1.074.$

PARTITIONS-----------------------

\subsection{Near-best spline iQIs}

A similar study can be done for integral spline QIs. We refer to
\cite{BISS1}\cite{IP} and we only
give  an example given in these papers. Setting $a=(a_0,a_1,\ldots,a_n)$,
$n\ge m$ and
$M_i(x)=M_{2m}(x-i)$, we consider
$
Q_af=\sum_{i\in\ZZ} \Lambda f(i) M_i
$
with coefficient functionals
$$
\Lambda f(i)=a_0\langle f,M_i\rangle+\sum_{j=1}^{n}a_j\left(\langle
f,M_{i-j}\rangle+\langle
f,M_{i+j}\rangle\right).
$$
As in section 4.2, we have  $\Vert Q_a\Vert_{\infty}\le \nu(a)$
and we say that $Q^*=Q_{a^*}$ is a {\sl near best iQI} if
$
\nu(a^*)=\min\{\nu(a); a\in V_r\}.
$
There is existence, but in general not unicity, of solutions.

\noindent
{\sl Example:} cubic splines (see \cite{IP}).
There is a unique optimal solution for $n\ge 2$:
$$
a_0^*=1+\frac{2}{3n^2},\;\; a_n^*=-\frac{1}{3n^2},\;\; a_j^*=0 \;\; for
\;\; 1\le j\le n-1
$$
Moreover, for all $n\ge 4$, $\Vert Q^* \Vert_{\infty}\le 1+\frac{4}{3n^2}.$
Here are the first values of $\Vert Q^* \Vert_{\infty}$ and $\nu(a^*);$
$n=1: 1.5278\;\; \&\;\; 2.333;$
$n=2: 1.2778\;\;  \&\;\; 1.333;$
$n=3: 1.1481\;\;  \&\;\; 1.1482.$

-----------------

\section{Bivariate spline dQIs on uniform partitions}

\subsection{A general construction of dQIs}

Let $\varphi$ be any kind of bivariate B-spline on one of the two classical
three- or
four-directional meshes of the plane (e.g. box-splines, see
\cite{BHS},\cite{dBHR},\cite{Chu}).
Let $\Sigma=supp(\varphi)$ and $\Sigma^*=\Sigma\cap{\ZZ}^2$.
Let $a$ be the hexagonal (or lozenge=rhombus) sequence formed by the values
$\{\varphi(i),i\in\Sigma^*\}$. The associated central difference operator
$\mathcal{D}$ is an isomorphism of $\PP(\varphi)$, the maximal subspace of
"complete " polynomials in the space of splines $\mathcal{S}(\varphi)$
generated by the integer translates of the B-spline $\varphi$ (see \cite{dBHR},
\cite{Sab9},\cite{Sab11},\cite{Sab12}).
Computing the expansion of $a$ in some basis of the space of
hexagonal (or lozenge) sequences amounts to expand $\mathcal{D}$
in some basis of central difference operators.
Then, computing the formal inverse ${\mathcal{D}}^{-1}$ allows to define the dQI
$$
{\mathcal{Q}}f=\sum_{k\in {\ZZ}^2}{\mathcal{D}}^{-1}f(k)\varphi(\cdot-k)
$$
which is exact on ${\PP} (\varphi)$. Let us now give two examples which are
detailed in
\cite{IP}.

\subsection{Near-best spline dQIs on a three direction mesh}

{\sl Example:} let $\varphi$ be the $C^2$ quartic box-spline. Let $H_s$ be
the regular hexagon
with edges of length $s\ge 1$, centered at the origin (here $\Sigma=H_2$)
and let
$H_s^*=H_s\cap {\ZZ}^2$. The near-best dQIs have coefficient functionals
with supports
consisting of the center and the 6 vertices of  $H_s^*, s\ge 1$. The
coefficients of values of
$f$ at those points are respectively
$1+\frac{1}{2s^2}$ and $-\frac{1}{12s^2}$, therefore the infinite norm of
the optimal dQIs
$Q_s^*$ is bounded above by $\nu_s^*=1+\frac{1}{s^2}$. Here are the first
values of
$\Vert Q^* \Vert_{\infty}$ and $\nu_s^*;$
$n=1: 1.34028\;\; \&\;\; 2 ;$
$n=2: 1.22917\;\;  \&\;\; 1.25;$
$n=3: 1.10185\;\;  \&\;\; 1.111.$

\subsection{Near-best spline dQIs on a four direction mesh}

{\sl Example:} let $\varphi$ be the $C^1$ quadratic box-spline. Let
$\Lambda_s$ be the lozenge
(rhombus) with edges of length $s\ge 1$, centered at the origin, and let
$\Lambda_s^*=\Lambda_s\cap {\ZZ}^2$. The near-best dQIs have coefficient
functionals with supports
consisting of the center and the 4 vertices of  $\Lambda_s^*, s\ge 1$. The
coefficients of values
of $f$ at those points are respectively
$1+\frac{1}{2s^2}$ and $-\frac{1}{8s^2}$, therefore the infinite norm of
the optimal dQIs
$Q_s^*$ is bounded above by $\nu_s^*=1+\frac{1}{s^2}$. Here are the first
values of
$\Vert Q^* \Vert_{\infty}$ and $\nu_s^*;$
$n=1: 1.5\;\; \&\;\; 2 ;$
$n=2: 1.25\;\;  \&\;\; 1.25;$
$n=3: 1.111\;\;  \&\;\; 1.111.$

--------------------

\section{Univariate spline QIs on non uniform partitions}

\subsection{Uniformly bounded dQIs}

Let us only give an example: we start from a family of DQIs of degree $m$
which are {\sl exact on}
${\PP}_2$.
$$
Q_2f=\sum_{j\in J} \lambda_j^{(2)}(f) B_j,\;\;
\lambda_j^{(2)}(f)=f(\theta_j)-\frac12 (\theta_j^2-\theta_j^{(2)})
D^2f(\theta_j).
$$
We recall the expansion \cite{MS1}\cite{MS2}
$$
A_j^{(2)}=\theta_j^2-\theta_j^{(2)}=\frac{1}{(m-1)^2 (m-2)}\sum_{(r,s)\in
E_m^2,r\neq
s}(t_{j+r}-t_{j+s })^2>0.
$$
On the other hand, $\frac12 D^2f(\theta_j)$ can be replaced on the space
${\PP}_2$ by the
second order divided difference $[\theta_{j-1}, \theta_j,\theta_{j+1}]f$,
therefore the dQI defined by
$$
Q_2^*f=\sum_{j\in J} \mu_j^{(2)}(f) B_j,\;\;
\mu_j^{(2)}(f)=f(\theta_j)-A_j^{(2)} [\theta_{j-1}, \theta_j,\theta_{j+1}]f,
$$
is also exact on ${\PP}_2$. Moreover, one can write
$$
\mu_i^{(2)}(f)=a_if_{i-1}+b_i f_i+c_i f_{i+1}
$$
with
$
a_i=-A_i^{(2)}/\Delta \theta_{i-1}(\Delta \theta_{i-1}+\Delta \theta_{i}),
\;\;
c_i=-A_i^{(2)}/\Delta \theta_{i}(\Delta \theta_{i-1}+\Delta \theta_{i}),$
 and \\
$b_i=1+ A_i^{(2)}/\Delta \theta_{i-1}\Delta \theta_{i}$,
So, according to the introduction
$$
\Vert Q_2^* \Vert_{\infty}\le \max_{i\in J} (\vert a_i \vert+\vert b_i
\vert +\vert c_i \vert)
\le 1+2\max_{i\in J}\frac{A_i^{(2)}}{\Delta \theta_{i-1}\Delta \theta_{i}}.
$$
The following theorem \cite{BISS3} extends a result given for quadratic
splines in
\cite{BISS3}\cite{Sab13}\cite{Sab15}.

\smallskip
\noindent

\proclaim Theorem 1.
For any degree $m$, the dQIs $Q_2^*$ are UB. More specifically, for all
partitions of $I$:
$$
\Vert Q_2^* \Vert_{\infty}\le [\frac12 (m+4)]
$$
---------------------
\subsection{Uniformly bounded iQIs}

General types of integral QIs  are studied in \cite{Cz}\cite{Sab3}\cite{Sab8}.
Here, we have chosen to study a family of QIs that we call Goodman-Sharma
type iQIs, as they
first appear  in \cite{GS1}. They seem simpler and more interesting  than
those we have studied in
\cite{Sab8}. The simpler GS-type IQI can be written as follows
$$
G_1f=f(t_0)B_0+\sum_{i=1}^{n+m-2}\tilde\mu_i(f)B_i+f(t_n)B_{n+m-1},
$$
where the integral coefficient functionals are defined by
$$
\tilde\mu_i(f)=\int_0^1 \tilde M_{i-1}(t)f(t)dt,
$$
$\tilde M_{i-1}(t)$ being the B-spline of degree $m-2$ with support
$\tilde \Sigma_{i-1}=[t_{i-m+1},t_i]$, normalized by $\tilde
\mu_{i}^{(0)}=\tilde \mu_{i}(e_0)=
\int_0^1\tilde M_{i-1}(t)=1$. It is easy to verify that $G_1$ is exact on
${\PP}_1$ and that $\Vert
G_1\Vert_{\infty}=1$. We shall study the family of GS-type iQIs defined by
$$
G_2f=f(t_0)B_0+\sum_{i=1}^{n+m-2}[a_i\tilde\mu_{i-1}(f)+b_i\tilde\mu_i(f)+c_
i\tilde\mu_{i+1}(f)]B_i
+f(t_n)B_{n+m-1},
$$
which are {\sl exact} on ${\PP}_2$. The three constraints
$G_2e_k=e_k,\;\;k=0,1,2$, lead to
the  following system of equations, for $1\le i\le n+m-2$:
$$
a_i+b_i+c_i=1,\quad
\theta_{i-1}a_i+\theta_{i}b_i+\theta_{i+1}c_i=\theta_{i},\quad
\tilde\mu_{i-1}^{(2)}a_i+\tilde\mu_{i}^{(2)}b_i+\tilde\mu_{i+1}^{(2)}c_i=\theta_i^{(2)}.
$$
This is a consequence of the following facts
$$
\tilde\mu_i(e_1)=\int_0^1 t\tilde
M_{i-1}(t)dt=\frac{1}{m}\sum_{s=1}^{m}t_{i-m+s}=\theta_i,
$$
$$
\tilde\mu_i^{(2)}=\mu_i(e_2)=\int_0^1 t^2\tilde
M_{i-1}(t)dt=\frac{2}{m(m+1)}\tilde s_2(T_i)
$$
$$
=\frac{2}{m(m+1)}\sum_{1\le r\le s\le m}t_{i-m+r}t_{i-m+s}
$$
\proclaim Theorem 2.
For any degree $m$, the iQIs $G_2$ are UB. More specifically, for all
partitions of $I$:
$$
\Vert G_2 \Vert_{\infty}\le 5
$$

The detailed proof will be given in \cite{Sab18}.

---------------------
\subsection{Near-best dQIs}

Let us consider the family of dQIs of degree $m$ defined, for the sake of
simplicity, on
$I=\RR$ endowed with an
arbitrary non-uniform increasing sequence of knots $T=\{t_i; i\in {\ZZ}\}$, by
$$
Qf=Q_{p,q}f=\sum_{i\in\ZZ} \mu_i(f) B_i.
$$
Their coefficient functionals depend on $2p+1$ parameters, with $p\ge m$:
$$
\mu_i(f)=\sum_{s=-p}^{p}\lambda_i(s)f(\theta_{i+s}),
$$
and they are exact on the space ${\PP}_q$, where $q\le \min(m,2p)$.
The latter condition is equivalent to $Q e_r=e_r$ for all monomials of
degrees $0\le r\le q$.
It implies that for all indices $i$, the parameters $\lambda_i(s)$ satisfy
the system of $q+1$
linear equations:
$$
\sum_{s=-p}^{p}\lambda_i(s)\theta_{i+s}^r=\theta_i^{(r)},\quad 0\le r\le q.
$$
The matrix $V_i\in {\RR}^{(q+1)\times (2p+1)}$ of this system, with
coefficients
$V_i(r,s)=\theta_{i+s}^r$, is a Vandermonde
matrix of maximal rank $q+1$, therefore there are $2p-q$ {\sl free parameters}.
Denoting $b_i\in {\RR}^{q+1}$ the vector in the right hand side, with
components
$b_i(r)=\theta_i^{(r)},\quad 0\le r\le q$, we consider the sequence of
minimization problems,
for $i\in \ZZ$:
$$
\min \Vert \lambda_i \Vert_1,\quad V_i \lambda_i=b_i.
$$
We have seen in the introduction that $\nu_1^*(Q)=\max_{i\in\ZZ}\min\Vert
\lambda_i \Vert_1$ is an upper bound
of $\Vert Q_q \Vert_{\infty}$ which is easier to evaluate than the true
norm of the dQI.

\proclaim Theorem 3. The above minimization problems have always solutions,
which, in general,
are non unique.

\noindent
The {\sl objective function being convex} and {\sl the domains being affine
subspaces},
these classical optimization problems have always solutions, in general non
unique.

Example of optimal dQIs are given in \cite{BIS}\cite{BISS3}\cite{IP}.

---------------------

\section{Bivariate quadratic spline dQIs on non uniform criss-cross
triangulations}

At the author's knowledge, the only bivariate box-splines which have been
extended to non uniform
partitions of the plane are $C^1$-quadratic box-splines on criss-cross
triangulations
\cite{CSW}\cite{Sab2}. Recently, we have constructed a set of B-splines
generating the space
of quadratic splines on a rectangular domain and we have defined a
discrete quasi-interpolant which is exact on ${\PP}_2$ and uniformly
bounded independently of
the partition \cite{Sab14}-\cite{Sab16}.

---------------------
\section{Abbreviations for publishers  and journals}

{\bf Publishers:} AP=Academic Press, New-York ; BAS=Bulgarian Academy of
Science, Sofia;
BV=Birkh\"auser-Verlag, Basel;
CUP=Cambridge University Press;\\
JWS=John Wiley \& Sons, New-York; K=Kluwer, Dordrecht; NH= North-Holland,
Amsterdam;
NP=Nashboro Press, Brentwood;
SV=Springer-Verlag, Berlin;\\ SIAM=Society for Industrial and Applied
Mathematics,
Philadelphia;\\ VUP=Vanderbilt University Press, Nashville.\\
\noindent
{\bf Journals:} AiCM=Advances in Comput. Mathematics; ATA=Approximation
Theory and its
Applications (now Analysis in Theory and Applications); \\CAGD=Computer
Aided Geometric
Design; JAT=Journal of Approximation Theory;  JCAM=Journal of Computational
and Applied
Mathematics.\\
\noindent
{\bf Proceedings:} AT2={\sl Approximation Theory II}, G.G.
Lorentz, C.K. Chui, L.L. Schumaker (eds), AP 1976;
AT4 \& AT5={Approximation Theory IV \& V}, C.K. Chui, L.L. Schumaker, J.D.
Ward (eds), AP 1983
and 1986; CMSB=Colloquia Mathematica Soc. Janos Bolyai;
CS02=Curve and Surface Fitting (St Malo
2002), A. Cohen, J.L. Merrien and L.L. Schumaker (eds), NP 2003.\\
\noindent
{\bf Preprints:} PI=Pr\'epublications IRMAR, Inst. de Recherche Math. de Rennes.



%


%


\begin{flushleft}

Centre de Math\'ematiques\\

INSA de Rennes\\

20 \ Avenue des Buttes de Co\"esmes, CS 14315\\

35043-Rennes cedex, France\\

{\em Email address:} psablonn@insa-rennes.fr

\end{flushleft}

\end{document}